\documentclass{amsart}

\usepackage{graphicx}
\usepackage{times}

\def \diam {\text{diam}}

\newtheorem{thm}{Theorem}[section]

\newtheorem{remark}[thm]{Remark}

\begin{document}

\title{The Radio Number of Gear Graphs}
\date{\today}

\author{Christina Fernandez}
\address{Department of Mathematics, CSU Channel Islands} \email{ma.christy@teacher.com}
\author{Am\'erica Flores}
\address{Department of Mathematics, CSU Channel Islands} \email{america.flores600@dolphin.csuci.edu}
\author{Maggy Tomova}
\address{Department of Mathematics, University of Iowa} \email{mtomova@math.uiowa.edu}
\author{Cindy Wyels}
\address{Department of Mathematics, CSU Channel Islands} \email{cindy.wyels@csuci.edu}


\begin{abstract}

Let $d(u,v)$ denote the distance between two distinct vertices of a connected graph $G$, and $\diam(G)$ be the diameter of $G$. A \emph{radio labeling} $c$ of $G$ is an assignment of positive integers to the vertices of $G$ satisfying $d(u,v)+|c(u)-c(v)|\geq \diam(G) + 1.$ The maximum integer in the range of the labeling is its span. The \emph{radio number} of $G$, $rn(G)$, is the minimum possible span. The family of gear graphs of order $n$, $G_n$, consists of planar graphs with $2n+1$ vertices and $3n$ edges.  We prove that the radio number of the $n$-gear is $4n+2$.

\vspace{.1in}\noindent \textbf{2000 AMS Subject Classification:} 05C78 (05C15)

\end{abstract}

\maketitle

\section{Introduction}
Radio labeling of graphs is motivated by restrictions inherent in assigning channel frequencies for radio transmitters \cite{CEHZ}.  To
avoid interference, transmitters that are geographically close must be assigned channels with large frequency differences; transmitters that are further apart may receive channels with relatively close frequencies.
The general situation is modeled by identifying transmitters with the vertices of a graph.  We assign positive integers to the vertices of the graph subject to a restriction concerning the distances between vertices;  the goal is to minimize the largest integer used.

We will consider simple connected graphs $G=(V(G), E(G))$.  We write $d(u,v)$ for the distance between vertices $u$ and $v$, and use $\diam(G)$ to indicate the diameter of $G$.  A \emph{radio labeling } is a one-to-one mapping $c: V(G)\rightarrow \mathbf {Z_+}$ satisfying the condition
\begin{equation}  \label{E:first}
d(u,v)+|c(u) - c(v)|\geq \diam(G) + 1
\end{equation}
for every $u,v\in V(G)$.  The \emph{span }of a labeling $c$ is the maximum integer that $c$
maps to a vertex of graph $G$. The \emph{radio number} of $G$, $rn(G)$,
is the lowest span taken over all radio labelings of the graph $G$\footnote{Some authors allow 0 in the range of a radio labeling (e.g. \cite{Liu1});  their radio numbers will be one less than what we define.}.  We will refer to Inequality \ref{E:first} as {\it the radio condition}.  Note that this condition necessitates the use of distinct integers, thus $rn(G) \geq |V(G)|$ for all graphs $G$.
Radio labelings are sometimes referred to as multi-distance labelings (e.g. \cite{Liu1}), and they are equivalent to $k$-labelings for $k=\diam(G)$.

In this introduction we will briefly note the radio numbers of some common families of graphs:  complete graphs, stars, and wheels.  The main section is devoted to establishing the radio number of gear graphs, which are extensions of wheel graphs.

%*********************************************Complete graphs**************************************

\begin{thm} The radio number of the complete graph on $n$ vertices is $n$, i.e. $rn(K_n)=n$.
\end{thm}
\begin{proof}
Since $|V(K_n)| = n$, we have $rn(K_n)\geq n$.  As $\diam(K_n)=1$, we may satisfy the radio condition while labeling the vertices with consecutive integers.
\end {proof}

\begin{figure}[tbh]
\centering
\includegraphics[scale=0.5]{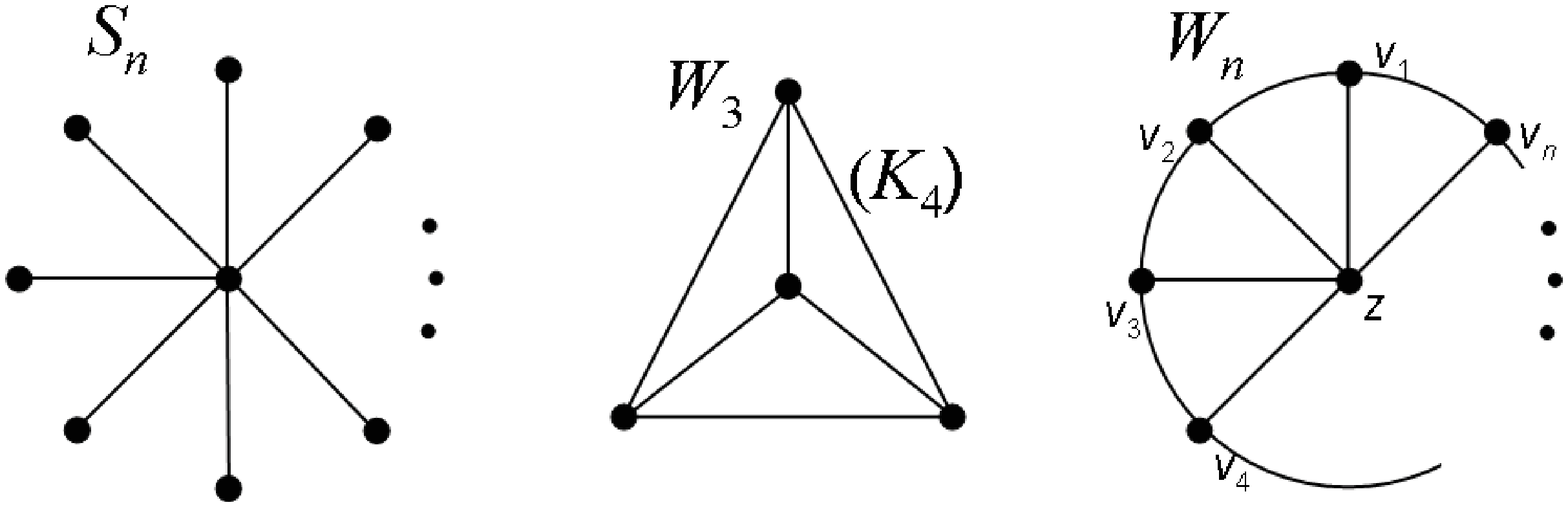}
\caption{} \label{fig:starwheels}
\end{figure}

%**********************************Star graphs*********************************************
The star graph $S_n$ ($n \geq 2$) is a tree on $n+1$ vertices.  One vertex (the ``center") is adjacent to every vertex; all other vertices have degree 1.  That is, $S_n = K_{1,n}$.  (See Figure \ref{fig:starwheels}.)

\begin{thm} $rn(S_n)= n+2$. \label{thm:star}
\end{thm}
\begin{proof}
Note that $\diam(S_n) = 2$.  This together with the fact that the center vertex is adjacent to every other vertex implies we may not use consecutive integers to label the center and another vertex.  Since $|V(S_n)| = n+1$, we see $rn(S_n) \geq n+2$.  Assigning 1 to the center and consecutive integers beginning with 3 to the other vertices produces a radio labeling with span $n+2$, so $rn(S_n)=n+2$.
\end{proof}

\begin{remark}  $rn(K_{m,n})=m+n+1$.
\end{remark}
\begin{proof}
Using similar logic as in the proof of Theorem \ref{thm:star}, we see that $rn(K_{m,n})\geq m+n+1$.  This lower bound may be realized by the span of a radio labeling:  use the label set $\{1, 2, \dots , m\}$ on the partition of size $m$, and the label set $\{m+2, m+3, \dots , m+n+1\}$ on the other partition.
\end{proof}
%***********************************Wheels*************************************************
The \textit{wheel graph} $W_n$ ($n \geq 3$) consists of an $n$-cycle together with a center vertex that is adjacent to all $n$ vertices of the cycle.  (See Figure 1.)  $W_3 = K_4$, so $rn(W_3)=4$.  As $\diam(W_4)=2$, adjacent vertices must have non-consecutive labels.  There are only two mutually-exclusive pairs of non-adjacent vertices, so at most two pairs of consecutive numbers may be used to label the five vertices.  This gives a lower bound for $rn(W_4)$ of 7;  a labeling assigning 1 to the center vertex and 3, 6, 4, and 7 (sequentially) to the vertices on the cycle realizes this bound and is a radio labeling.  Thus $rn(W_4)=7$.  The radio number of all larger wheels is given in Theorem \ref{thm:wheel}.

\begin{thm} \label{thm:wheel}
$rn(W_n)= n+2$ for $n \geq 5$.
\end{thm}
\begin{proof}
For $n\geq 4$, $\diam(W_n)=2$.  As with the star, we may not use consecutive integers to label the center and another vertex, so $rn(W_n) \geq n+2$.  To define a radio labeling, we name the center vertex $z$, and identify the remaining vertices as $\{v_1, v_2, \dots v_n\}$, labeled sequentially around the cycle (as shown in Figure \ref{fig:starwheels}).  Our labeling $c$ assigns $z$ the label 1.  Vertices $v_1, v_2, \dots v_{\lceil \frac{n}{2}\rceil}$ are labeled with consecutive odd numbers, beginning with 3.  Vertices $v_{\lceil\frac{n}{2}\rceil+1}, \dots v_n$ are assigned consecutive even integers beginning with 4.

We must show that the radio condition is satisfied for all pairs of vertices.  Consider the pair $(z,v_i)$, noting that $d(z,v_i) = 1$.  We have $c(z)=1$ and $c(v_i) \geq 3$.  So
$$d(z,v_i) + |c(z)-c(v_i)| \geq 1+|3-1|\geq 3 = \diam(W_n)+1.$$
If $v_i$ and $v_j$ ($i \neq j$) are adjacent, their labels differ by at least 2.  Otherwise, $d(v_i,v_j)\geq 2$, and $|c(v_i)-c(v_j)| \geq 1$.  In either case, the radio condition is satisfied.

As the range of $c$ is the label set $\{1, 3, 4, \dots , n+2\}$, this establishes $rn(W_n)=n+2$.
\end{proof}

%*********************************Gear graphs***********************************
\section{The Radio Number of Gear Graphs}

The proof of Theorem \ref{thm:wheel} models the strategy we will use to establish the radio number of gear graphs.  The lower bound is found by examining the minimum necessary differences between labels.  To determine an upper bound, we provide a specific radio labeling.  As the labeling provided has span equal to the lower bound, their common value is the radio number.

Gear graphs are extensions of wheel graphs.  We may describe the $n$-gear, $G_n$, as $W_n$ with an additional vertex inserted between each pair of vertices on the cycle.  Equivalently, the $n$-gear consists of a cycle on $2n$ vertices, with every other vertex on the cycle adjacent to an additional $2n+1^{\text{st}}$ vertex (the center).  Gear graphs have $3n$ edges.  For $n\geq 4$, $\diam(G_n)=4$.

\begin{figure}[tbh]
\centering
\includegraphics[scale=0.5]{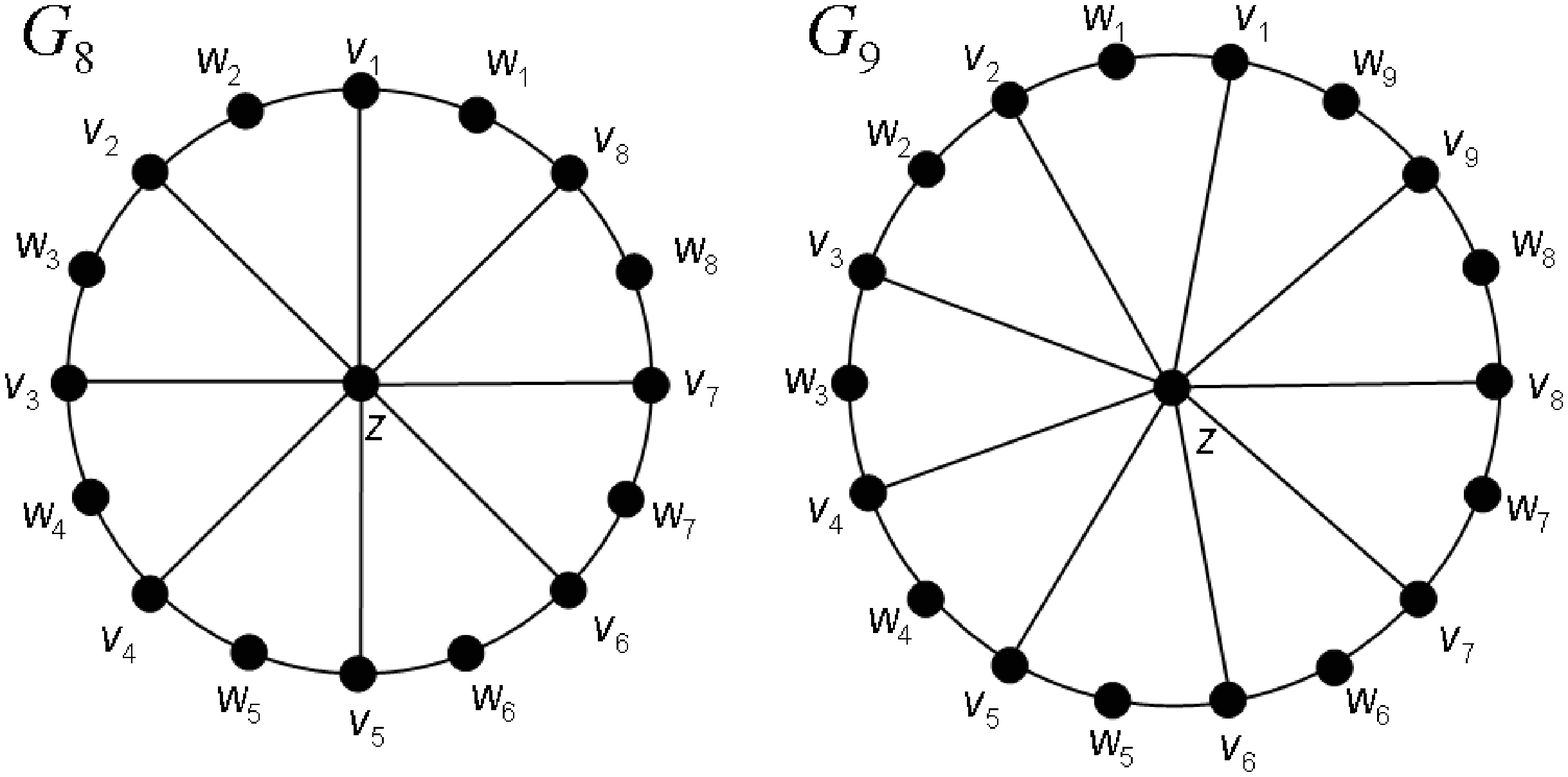}
\caption{} \label{fig:standard}
\end{figure}

%*******************Standard labeling for the gear graphs*************************************
 \vspace{.1 in} \noindent{\it The Standard Labeling of the $n$-Gear}
\newline\noindent
To establish the radio number of the $n$-gear we will refer to a labeling of the vertices of the $n$-gear that distinguishes the vertices by their characteristics.  The center vertex is labeled $z$, and the vertices adjacent to the center are labeled sequentially $\{v_1, v_2, ...,v_n\}$.  The vertices not adjacent to the center are labeled sequentially $\{w_1, w_2,...w_n\}$, using the same orientation chosen for the $v$'s.  If $n$ is odd then we specify that $w_1$ is adjacent to $v_1$ and $v_2$, otherwise $w_1$ is adjacent to $v_1$ and $v_n$.  The standard labelings of the 8-gear and of the 9-gear are shown in Figure \ref{fig:standard}.

%**************************Lower bound*******************************************************
\begin{thm} \label{thm:lowerbound}
For $n\geq 4, rn(G_n)\geq 4n+2$.
\end{thm}
\begin{proof}
Assume $n \geq 4$.  This gives $\diam(G_n)=4$, so any radio labeling $c$ of $G_n$ must satisfy the radio condition
$$d(u,v)+|c(u) - c(v)|\geq 5$$
for all distinct $u,v \in V(G)$.
We count the number of values needed for labels and add the minimum number of {\it forbidden} values -- those values precluded by use of a particular value as a label.  For instance, if we label the center $a$ (i.e. $c(z) = a$), then, as $d(z,r) \leq 2$ for all vertices $r \neq z$, the values $a-2$, $a-1$, $a+1$, and $a+2$ are forbidden.  Similarly, as $d(v_i,r) \leq 3$ for all $v_i$ and for any $r \neq v_i$, one value is forbidden above and below any label $c(v_i)$.  However, as $d(w_i,r) = 4$ for some vertex $r$, it is possible to use consecutive labels on $w_i$ and $r$.  (I.e.~there are no forbidden values associated with the vertices $\{w_1, \dots w_n\}$.)

Note that the number of forbidden values is symmetric below and above any label used for a particular vertex.  Thus we find the minimum number of forbidden values by using the lowest and highest-valued labels on the center vertex and on one of $\{v_1, \dots v_n\}$.  Assume without loss of generality that $c(z)=1$ and $c(v_n)$ is the span of $c$.  Associated with the center vertex are then two forbidden values (2 and 3), with $v_n$ is one forbidden value (the span less one), and with the other $v_i$ are two forbidden values each.  This gives a total of $2 + 1 + 2(n-1) = 2n+1$ forbidden values.  Adding in the $2n+1$ values needed to label the $2n+1$ vertices provides a total of $4n+2$, hence $rn(G_n) \geq 4n+2$.

\end{proof}

%******************Upper Bound**************************************************
\begin {thm}For $n\geq 7, rn(G_n)\leq 4n+2$. \label{thm:upperbound}
\end {thm}
\begin{proof}
We provide a radio labeling $c$ of $G_n$ in two steps:  first we define a position function that renames the vertices of $G_n$ using the set $\{x_0,x_1, \dots x_{2n}\}$, then we specify the labels $c(x_i)$ so that $i < j$ if and only if $c(x_i)<c(x_j)$.  (This allows us to more easily show that $c$ is indeed a radio labeling.) Once it is established that $c$ is a radio labeling, the span of $c$ provides an upper bound for $rn(G_n)$. Throughout this proof, $n \geq 7$.

The position function $p:V(G_n) \rightarrow \{x_0,x_1, \dots x_{2n}\}$ is defined as follows.
For $n=2k+1$ we define
$$
\begin{array}{lll}
p(z)&=x_0, & \\
p(w_{2i-1})&=x_i &\text{ for } i = 1, \dots , k+1, \\
p(w_{2i})&=x_{k+1+i} &\text{ for } i=1, \dots k, \\
p(v_i) &= x_{n+i} &\text{ for } i = 1, \dots, n.
\end{array}
$$
When $n=2k$ the position function changes slightly in renaming the vertices $w_i$:
$$
\begin{array}{lll}
p(z)&=x_0, & \\
p(w_{2i-1})&=x_i &\text{ for } i = 1, \dots , k, \\
p(w_{2i})&=x_{k+i} &\text{ for } i=1, \dots k, \\
p(v_i) &= x_{n+i} &\text{ for } i = 1, \dots, n.
\end{array}
$$
In essence, the position function orders the vertices so that $\{x_0, x_1, \dots x_{2n}\}$ corresponds to $\{z, w_1, w_3, \dots ,w_n, w_2, w_4, \dots , w_{n-1}, v_1, v_2, \dots v_n\}$ when $n$ is odd and to 

\noindent $\{z, w_1, w_3, \dots ,w_{n-1}, w_2, w_4, \dots , w_n, v_1, v_2, \dots v_n\}$ when $n$ is even.  Figure \ref{fig:position} shows the renamed versions of the 8-gear and the 9-gear.

\begin{figure}[tbh]
\centering
\includegraphics[scale=0.5]{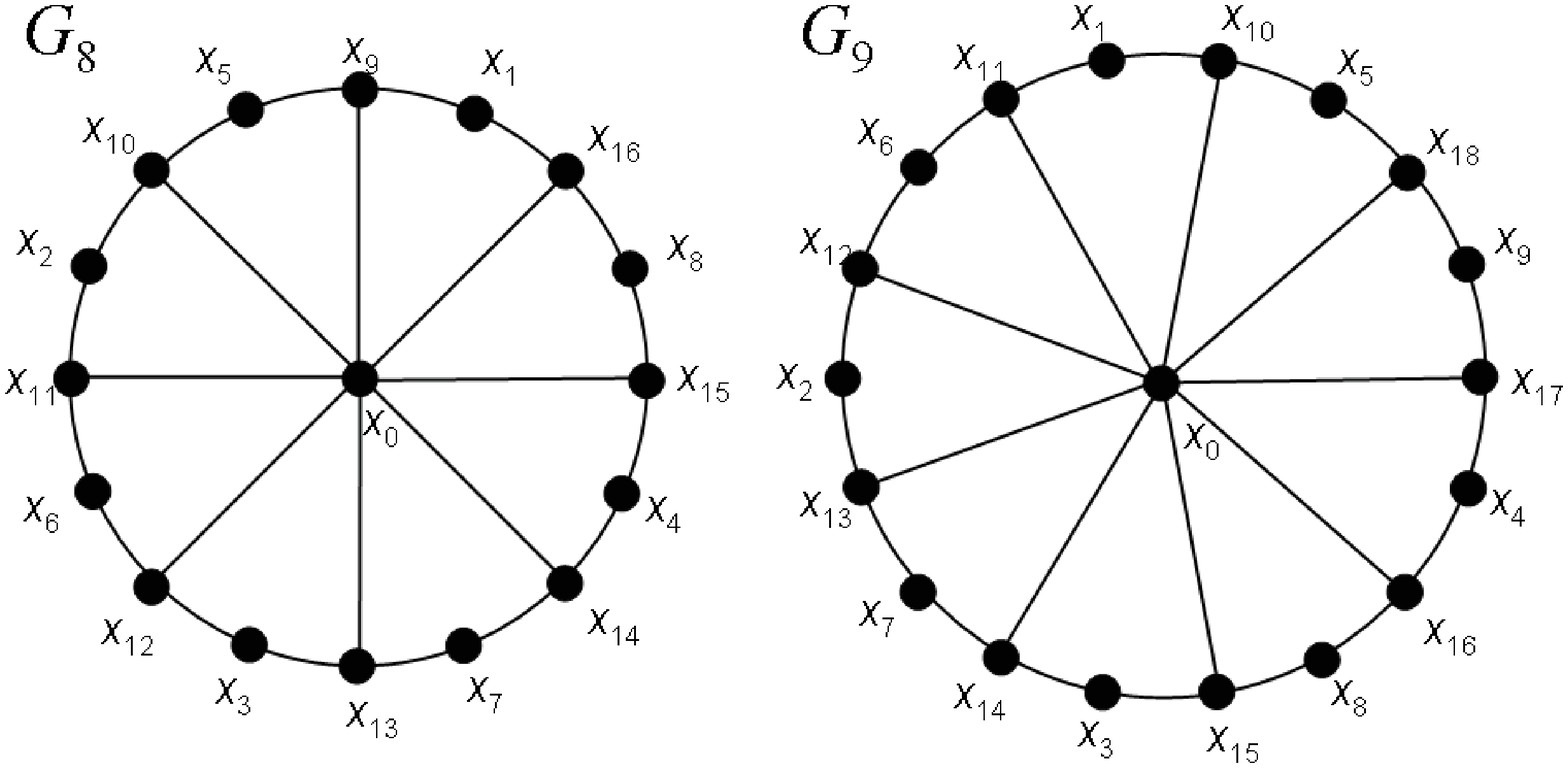}
\caption{} \label{fig:position}
\end{figure}

We are now ready to define our radio labeling $c:V(G)\to\mathbf {Z_+}$.
$$c(x_i)=\left\{
\begin{array}{ll}
1, & i=0, \\
3+i, & 1\leq i\leq n, \\
n+2+3(i-n), & n+1\leq i \leq 2n.
\end{array}
\right.$$
Figure \ref{fig:position} shows the labeling $c$ applied to the 8-gear and to the 9-gear.

\noindent\textbf{Claim:} The labeling $c$ is a valid radio labeling.
We must show that the radio condition,
$$d(u,v)+ |c(u)-c(v)|\geq \diam(G_n)+1=5,$$ holds for all pairs of vertices $(u,v)$ (where $u \neq v$).

\textbf{Case 1}:  Consider the pair $(z,r)$ (for any vertex $r \neq z$).
\newline
Recall $p(z)=x_0$.  As $c(x_i)\geq 5$ for any $i \geq 2$, we have $d(z,x_i)+|c(z)-c(x_i)|\geq 1+|1-5|\geq 5$ for all $i \geq 2$.  This leaves the pair $(z,x_1)$.  But $p^{-1}(x_1)=w_1$, so we calculate $d(z,w_1)+|c(z)-c(w_1)|=2+|1-4|=5$.

\textbf{Case 2}: Consider the pair $(w_j, w_k)$ (with $j \neq k$).
\newline
Recall $p(w_{2i-1})= x_i$ and note that $p(w_{2i})$ may be written as $x_{n-k+i}$ whether $n$ is even or odd.  We have $d(w_j,w_k)=2$ for the pairs $(w_{2i-1},w_{2i})$, $(w_{2i},w_{2i+1})$, and $(w_n,w_1)$.  These ``translate", respectively, to $(x_i,x_{n-k+i})$, $(x_{n-k+i},x_{i+1})$, and $(x_s,x_1)$, where $s=k+1$ when $n$ is odd and $s=2k$ when $n$ is even.  Examine the label difference for each pair:  $|c(x_i)-c(x_{n-k+i})| = n-k$, $|c(x_{n-k+i})-c(x_{i+1})|=n-k-1$, and $|c(x_s)-c(x_1)|$ is $k$ when $s=k+1$ ($n$ odd) and is $2k-1$ when $s=2k$ ($n$ even).  In all cases, using the fact that $n\geq 7$, we have that $|c(w_j)-c(w_k)| \geq 3$, so the radio condition is satisfied whenever $d(w_j,w_k)=2$.  Meanwhile, should $j$ and $k$ not be consecutive (mod $n$), we have $d(w_j,w_k)\geq 4$.  As it is always the case that $|c(w_j)-c(w_k)|\geq 1$, the radio condition is again satisfied.

\textbf{Case 3}: Examine the pair $(v_j,v_k)$ (with $j \neq k$).
\newline
Note that $d(v_j,v_k)=2$. Since $|c(v_j)-c(v_k)| = |c(x_{n+j})-c(x_{n+k})|\geq 3$ for all
$v_j,v_k$, the radio condition is satisfied.

\textbf{Case 4}: Finally, consider the pair $(v,w)$, where $v\in\{v_1,...,v_n\}$ and $w\in\{w_1,...,w_n\}$.
\newline
We have $c(v) \in \{n+5,n+8,...,4n+2\}$ and $c(w) \in \{4,5,6...,n+3\}$.  For all $v\neq v_1$, $|c(v)-c(w)|\geq (n+8)-(n+3)=5$.  Therefore the radio condition is satisfied when $v\neq v_1$.  Meanwhile, $|c(v_1)-c(w)|\geq (n+5)-(n+3)=2$.  If $w$ is distance three or greater from $v_1$, the radio condition holds.
If $d(v_1,w)<3$, then $w=x_1$ or $w=x_{\lfloor \frac n2 \rfloor +1}$.  Checking the radio condition for each we obtain $d(v_1,x_1)+|c(v_1)-c(x_1)| = 1 + |(n+5)-4| = n+2 \geq 5$, and $d(v_1,x_{\lfloor \frac n2 \rfloor +1})+|c(v_1)-c(x_{\lfloor \frac n2 \rfloor +1})| = 1 + |(n+5)-(3+{\lfloor \frac n2 \rfloor +1})| \geq \frac n2+3 \geq 5$.
%It remains to consider $|c(v_1)-c(w)|$ for $w$ adjacent to $v_1$:  our placement of $w_1$ adjacent to $v_1$ and $v_2$ when $n$ is odd and adjacent to $v_1$ and $v_n$ when $n$ is even ensures that $|c(v_1)-c(w)|\geq 5$ when $w$ is adjacent to $v_1$.

These four cases establish the claim that $c$ is a radio labeling of $G_n$.  Thus $rn(G_n) \leq \text{span}(c) = c(2n) = n+2+3(2n-n) = 4n+2$.
\end{proof}

Taken together, Theorem \ref{thm:lowerbound} and Theorem \ref{thm:upperbound} establish the main result of this paper:
\begin{thm}
The radio number of the $n$-gear is $4n+2$ when $n \geq 4$.
\end{thm}
\begin{proof}
Theorem \ref{thm:lowerbound} shows $rn(G_n)\geq 4n+2$ for $n \geq 4$, and Theorem \ref{thm:upperbound} shows $rn(G_n)\leq 4n+2$ for $n \geq 7$.  It remains only to show that $rn(G_n)\leq 4n+2$ for $n = 4, 5, 6$.  This is demonstrated by the radio labelings provided in Figure \ref{fig:smallgears}.  Note that the labels with values $a$ and $a+2$ are assigned to vertices at distance 3.
\end{proof}

For completeness' sake, we provide radio labelings of $G_n$ for $n = 2, \dots , 6$ in Figure \ref{fig:smallgears}.  The reader may verify that each radio labeling provided uses the minimum possible span, and thus exhibits the radio number of each small gear.

\begin{figure}[tbh]
\centering
\includegraphics[scale=0.5]{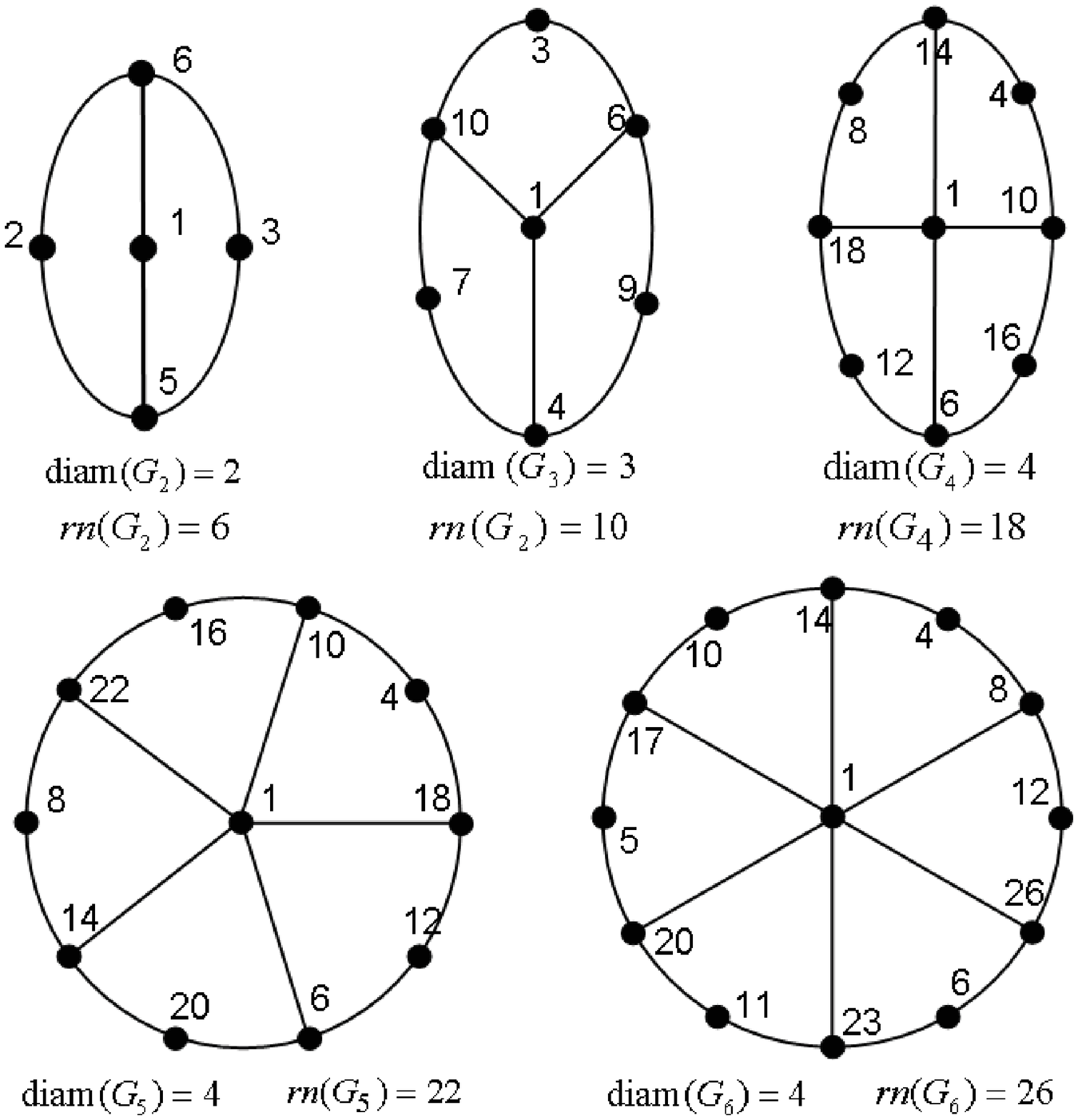}
\caption{} \label{fig:smallgears1}
\end{figure}

%\begin{remark}
%An alternative labeling with span $4n+2$ is easy to define for all $n$-gears with $n \not\equiv 0 \mod 3$ and $n \geq 4$.
%\end{remark}
%As before, we call the center vertex $x_0$ and define $c(x_0)=1$.  Choose any degree-2 vertex $x_1$ on the cycle and assign $c(x_1)=4$.  Traveling in one direction around the cycle, assign a vertex $x_2$ at distance 3 from $x_1$ the label 6.  Continue:  $x_i$ is the vertex at distance 3 from $x_{i-1}$ when traveling in the chosen direction, and $c(x_i)=2i+2$.  As $n \not\equiv 0 \mod 3$, all vertices receive a distinct label.  We see that $d(x_0,x_1)+|c(x_0)-c(x_1)| = 2 + 3$, and $|c(x_0)-c(x_i)| \geq 5$ for all $i \geq 2$.  Furthermore, for $i, j \not= 0$ ($i \not= j$), $d(x_i,x_j)=3$ implies $|c(x_i)-c(x_j)|=2$;  when $d(x_i,x_j)\not= 3$, we have $|c(x_i)-c(x_j)|\geq 4$.  The last vertex labeled is $x_{2n}$;  its label value is $c(x_{2n})=4n+2$.
%This labeling is demonstrated for $G_4$ and $G_5$ in Figure \ref{fig:smallgears}.

\section{Acknowledgements}
This research was carried out under the auspices of an MAA (SUMMA) Research Experience for Undergraduates program funded by NSA, NSF, and Moody's, and hosted at CSU Channel Islands during Summer, 2006.  We are grateful to all for the opportunities provided.
%*******************************Bibliography begins here.***************************************

\end{document}